# *Key-avoidance for alternating sign matrices*

Mathilde Bouvel[1] Rebecca Smith[2] Jessica Striker[3*]

[1] *Université de Lorraine, CNRS, Inria, LORIA, F-54000 Nancy, France*
[2] *Department of Mathematics, SUNY Brockport, Brockport, NY, USA*
[3] *Department of Mathematics, North Dakota State University, Fargo, ND, USA*

We initiate a systematic study of key-avoidance on alternating sign matrices (ASMs) defined via pattern-avoidance on an associated permutation called the *key* of an ASM. We enumerate alternating sign matrices whose key avoids a given set of permutation patterns in several instances. We show that ASMs whose key avoids 231 are permutations, thus any known enumeration for a set of permutation patterns including 231 extends to ASMs. We furthermore enumerate by the Catalan numbers ASMs whose key avoids both 312 and 321. We also show ASMs whose key avoids 312 are in bijection with the gapless monotone triangles of [Ayyer, Cori, Gouyou-Beauchamps 2011]. Thus key-avoidance generalizes the notion of 312-avoidance studied there. Finally, we enumerate ASMs with a given key avoiding 312 and 321 using a connection to Schubert polynomials, thereby deriving an interesting Catalan identity.

**Keywords:** alternating sign matrix, permutation pattern, pattern avoidance, ASM key, Catalan numbers

## 1 Introduction

Alternating sign matrices were introduced and investigated by Mills, Robbins, and Rumsey [12, 13], who conjectured a lovely enumeration formula. The search for and eventual completion of proofs of this enumeration [21, 8] motivated further investigation of their combinatorial, algebraic, and geometric properties and revealed a connection to statistical physics. Alternating sign matrices are in bijection with many interesting objects (see e.g. [14, 22]), including certain configurations of the square ice / six-vertex model, Gog triangles and Gog words, and are equinumerous with several other sets of combinatorial objects for which no explicit bijection is known.

In this paper, we initiate a study of key-avoidance of patterns in alternating sign matrices from an enumerative perspective. Our work is partly inspired by [5] in which this notion for a particular pattern was essential in finding a partial solution to one of these missing bijections. It is also motivated by a question left open in [2], where the authors define a notion of pattern-avoidance in Gog words, for a specific pattern only, and ask for a generalization of this notion to other patterns. Our work answers this question, on ASMs rather than on Gog words. (This is further discussed in Remark 3.9.)

**Definition 1.1.** An *alternating sign matrix* (or ASM) is a square matrix where

*Supported by a Simons Foundation gift MP-TSM-00002802 and NSF grant DMS-2247089.

1. Each entry is one of $0, 1, -1$.
2. Every row and column sum to 1.
3. The non-zero entries of each row and column alternate in sign.

ASMs can be thought of as a generalization of permutation matrices. Indeed, permutation matrices are exactly the alternating sign matrices with no $-1$ entries. The study of pattern-avoidance on permutations (or permutation matrices) has proved a very rich topic of combinatorics, with connections to sorting procedures, algorithms, formal logic, algebraic geometry or discrete probability theory, among others. We give the definition of permutation pattern-avoidance below, and then discuss extending the notion to alternating sign matrices.

**Definition 1.2.** A permutation $\sigma = \sigma(1) \dots \sigma(n)$ *contains* as a pattern the permutation $\pi = \pi(1) \dots \pi(k)$ if there exists indices $i_1, i_2, \dots, i_k$ such that, for all $a, b$, $\sigma(i_a) < \sigma(i_b)$ if and only if $\pi(a) < \pi(b)$. Equivalently, $\pi$ is a pattern of $\sigma$ when the permutation matrix of $\sigma$ contains as a submatrix the permutation matrix of $\pi$.
When $\sigma$ does not contain $\pi$, we say that $\sigma$ *avoids* $\pi$.

The purpose of this work is to investigate pattern-avoidance on ASMs. We will restrict our attention to patterns that are themselves permutations. There are (at least) two natural ways of extending the notion of pattern-avoidance from permutations to ASMs.

- One is to simply treat $-1$s like 0s in the alternating sign matrix and then determine whether or not the pattern (which is a permutation matrix) appears as a submatrix of this modified version of the alternating sign matrix. This notion was studied by Johansson and Linusson [6]. We refer to it as *classical avoidance*. We have a paper with Justin Troyka and Eric Egge [3] in which we further investigate this notion.
- The other, which we focus on in the present work, consists in associating a permutation to any ASM, and to define that the ASM avoids a (permutation) pattern if and only if the associated permutation does. A natural way to associate a permutation to an ASM is through its *key*. The process of computing the key of an ASM by systematically removing its $-1$s has been defined by Lascoux [9] and studied by Aval [1], and we review it below (with our specific convention) in Definition 1.3. Lascoux's motivation was to give formulas in *Schubert calculus* for *Schubert* and *Grothendieck* polynomials as a sum over weights associated to ASMs whose key is $\sigma$ [9]. Weigandt [19] later connected this perspective to the *bumpless pipe dreams* of Lam, Lee, and Shimozono [10].

First, we clarify our convention for drawing the permutation matrix of a permutation. Note that this is *not* the most common one in the permutation patterns literature, but we found that it is more practical here. To draw the permutation matrix of $\sigma = \sigma(1) \dots \sigma(n)$, for each $i \in \{1, \dots, n\}$, put a 1 in the $i$-th row from the top and in the $\sigma(i)$-th column from the left. So, for instance, the permutation matrix of 231 (resp. 3124) is $\begin{pmatrix} 0 & 1 & 0 \\ 0 & 0 & 1 \\ 1 & 0 & 0 \end{pmatrix}$ (resp. $\begin{pmatrix} 0 & 0 & 1 & 0 \\ 1 & 0 & 0 & 0 \\ 0 & 1 & 0 & 0 \\ 0 & 0 & 0 & 1 \end{pmatrix}$). Consistently, when we refer to entry in position $(i, j)$ in an ASM, we mean the entry in the $i$-th row from the top and the $j$-th column from the left.

We now turn to defining the key of an ASM. As indicated above, the key of an ASM is a permutation matrix obtained by systematically removing its $-1$ entries, along a preferred direction. Most classically (see e.g. [1] or more recent works on the topic), the removal of the $-1$ entries pushes the 1s to the northwest. In the context of our work, it makes sense to consider a symmetry of this construction, pushing the 1s to the southwest. The process of computing the key is illustrated in Fig. 1.

**Definition 1.3.** The *southwest key* (denoted *SW key* or *key* for short) of an alternating sign matrix $A$ is computed by removing the $-1$s one by one as follows (we sometimes refer to this process as the *SW key process*).

1. A $-1$ entry is removable if there are no other $-1$ entries that are weakly southwest of it.
2. Following Aval [1], we define the neighboring 1s of a removable $-1$ as the $1$ entries weakly to its southwest, for which no other $1$ entry lies both weakly to their northeast and weakly to the southwest of the considered $-1$. (In other words, they are the $1$ entries weakly to the southwest of the considered $-1$ and closest to it.)
3. For each removable $-1$,
   (a) Consider the Ferrer's shape(i) that results using the $-1$ as the northeast corner, the next $1$ entry to the west of the removable $-1$ as the northwest corner, the next $1$ to the south of the removable $-1$ as the southeast corner and any neighboring 1s of the removable $-1$ that appear in this rectangle as (inner) corners in the Ferrers diagram.
   (b) Replace the south-most 1 from the Ferrer's diagram with a $0$. Moving east to west, for each subsequent 1, place a new $1$ in the row of the previously replaced 1 and column (say, $c$) of the current 1, then replace the old 1 in column $c$ with a $0$. Finally, replace the $-1$ with a $0$.
4. The resulting permutation matrix determines the key of the ASM we started from.

Note if there are no inner corners in Step 3(a), we say the removal is *simple*.

As shown in [1, Corollary 9], the order in which the $-1$s are removed does not affect the resulting key.

**Definition 1.4.** We say that an ASM $A$ *key-avoids* a permutation (or permutation matrix) $\pi$ when the key of $A$ avoids $\pi$ in the sense of permutation pattern-avoidance.

**Definition 1.5.** Let $\pi$ be a permutation.

- We denote by $\mathrm{pAv}(\pi)$ the set of permutations that avoid $\pi$.
- We denote by $\mathrm{Av}(\pi)$ the set of ASMs that classically avoid $\pi$.
- We denote by $\mathrm{Av}^{\text{key}}(\pi)$ the set of ASMs which key-avoid $\pi$.

Throughout the paper, an index $n$ (as in $\mathrm{pAv}_n(\pi)$) indicates that we restrict to objects of size $n$. We also extend these notations to sets of permutation patterns $S$, where avoidance of $S$ means avoidance of each of the patterns in $S$.

(i) In this context, the Ferrers shape we want to consider have their first row on the top and their first column on the *right*; in other words, they are reflected along a vertical axis with respect to the classical convention.

$$\begin{pmatrix} 0 & 0 & 1 & 0 & 0 \\ 0 & 1 & -1 & 1 & 0 \\ 1 & 0 & 0 & -1 & 1 \\ 0 & 0 & 1 & 0 & 0 \\ 0 & 0 & 0 & 1 & 0 \end{pmatrix} \to \begin{pmatrix} 0 & 0 & 1 & 0 & 0 \\ 0 & 1 & -1 & 1 & 0 \\ \boxed{1} & \boxed{0} & \boxed{0} & \boxed{-1} & 1 \\ \boxed{0} & \boxed{0} & \boxed{1} & \boxed{0} & 0 \\ 0 & 0 & \boxed{0} & \boxed{1} & 0 \end{pmatrix} \to \begin{pmatrix} 0 & 0 & 1 & 0 & 0 \\ 0 & 1 & -1 & 1 & 0 \\ \boxed{0} & \boxed{0} & \boxed{0} & \boxed{0} & 1 \\ \boxed{1} & \boxed{0} & \boxed{0} & \boxed{0} & 0 \\ 0 & 0 & \boxed{1} & \boxed{0} & 0 \end{pmatrix}$$

$$\to \begin{pmatrix} 0 & 0 & 1 & 0 & 0 \\ 0 & \boxed{1} & \boxed{-1} & 1 & 0 \\ 0 & \boxed{0} & \boxed{0} & 0 & 1 \\ 1 & \boxed{0} & \boxed{0} & 0 & 0 \\ 0 & \boxed{0} & \boxed{1} & 0 & 0 \end{pmatrix} \to \begin{pmatrix} 0 & 0 & 1 & 0 & 0 \\ 0 & \boxed{0} & \boxed{0} & 1 & 0 \\ 0 & \boxed{0} & \boxed{0} & 0 & 1 \\ 1 & \boxed{0} & \boxed{0} & 0 & 0 \\ 0 & \boxed{1} & \boxed{0} & 0 & 0 \end{pmatrix} = 34512$$

**Fig. 1:** Determining the key of an alternating sign matrix.

Tables 1 and 2 summarize our results of Sections 2–4. Most importantly, we show both $\mathrm{Av}_n^{\mathrm{key}}(231)$ and $\mathrm{Av}_n^{\mathrm{key}}(312, 321)$ are enumerated by the Catalan numbers, the equality $|\mathrm{Av}_n^{\mathrm{key}}(132)| = |\mathrm{Av}_n^{\mathrm{key}}(213)|$ holds, and $\mathrm{Av}_n^{\mathrm{key}}(312)$ is in bijection with *gapless monotone triangles* of order $n$, studied by Ayyer, Cori, and Gouyou-Beauchamps [2].

Tables 1 and 2 also list the numbers of ASMs of small sizes key-avoiding each pattern of size 3 and each pair of patterns of size 3. We computed in Table 3 such numbers for all patterns of size 4, without obtaining any known sequences in the OEIS [17], and include this data for completeness. All new sequences found in these tables were computed using Sage [15].

| Excluded pattern | First terms of enumeration sequence | Known sequence |
|---|---|---|
| 123 | 1, 2, 6, 26, 156, 1276, 14109 | Open |
| 132 or 213 | 1, 2, 6, 27, 184, 1916, 30799 | Open, see Cor. 2.8 for equivalence |
| 231 | 1, 2, 5, 14, 42, 132, 429 | Catalan numbers, Prop. 2.3<br>OEIS A000108 |
| 312 | 1, 2, 6, 26, 162, 1450, 18626 | Gapless monotone triangles, Theorem 3.8<br>OEIS A180349 |
| 321 | 1, 2, 6, 25, 143, 1138, 12857 | Open |

**Tab. 1:** Enumeration sequences of ASMs key-avoiding permutation patterns of size 3.

In Section 5, we further study the structure of ASMs in $\mathrm{Av}^{\mathrm{key}}(312, 321)$, expressing their enumeration as a product of Catalan numbers, thereby deriving an interesting Catalan number identity.

The paper is organized as follows. In Section 2, we prove some initial enumerative results. In Section 3, we prove Theorem 3.8 showing that the gapless monotone triangles of [2] correspond to precisely to the 312-key-avoiding ASMs. In Section 4, we prove Theorem 4.8 that $\mathrm{Av}^{\mathrm{key}}(312, 321)$ is enumerated by the Catalan numbers. In Section 5, we decompose matrices in $\mathrm{Av}^{\mathrm{key}}(312, 321)$ into blocks and enumerate them using prior work related to Schubert polynomials and pipe dreams, deriving an unusual Catalan

| Excluded patterns | First terms of enumeration sequence | Known sequence |
|---|---|---|
| 123 & 132 or 123 & 213 | 1, 2, 5, 16, 63, 299, 1682 | Open<br>see Prop. 2.7 for equivalence |
| 123 & 231 | 1, 2, 4, 7, 11, 16, 22 | $\binom{n}{2}+1$, OEIS A000124<br>Corollary 2.6 |
| 123 & 312 | 1, 2, 5, 13, 35, 101, 325 | Open |
| 123 & 321 | 1, 2, 5, 9, 0, 0, 0, ... | Corollary 2.2 |
| 132 & 213 | 1, 2, 5, 18, 98, 823, 10866 | Open |
| 132 & 231 or 213 & 231 or 231 & 312 or 231 & 321 | 1, 2, 4, 8, 16, 32, 64 | $2^{n-1}$, OEIS A000079 (shifted), Corollary 2.6 |
| 132 & 312 or 213 & 312 | 1, 2, 5, 17, 80, 524, 4808 | Open<br>see Prop. 2.7 for equivalence |
| 132 & 321 or 213 & 321 | 1, 2, 5, 15, 59, 335, 2910 | Open<br>see Prop. 2.7 for equivalence |
| 312 & 321 | 1, 2, 5, 14, 42, 132, 429 | Catalan numbers, Theorem 4.8<br>OEIS A000108 |

**Tab. 2:** Enumeration sequences of ASMs key-avoiding pairs of patterns of size 3.

identity.

# 2 Initial enumeration results

We begin with some initial enumeration results, included for completeness (the main results of this paper are in Section 3 and following). From the famous Erdös-Szekeres [4] result on permutations avoiding both monotone increasing and monotone decreasing subsequences of the same size, we have the following proposition.

**Proposition 2.1.** *If a permutation* $\pi$ *avoids both* $123$ *and* $321$*, then* $\pi$ *must have size at most* $4$*.*

As the SW-key of an alternating sign matrix is a permutation, a direct consequence of this result is:

**Corollary 2.2.** *There is no ASM of size* $n > 4$ *that key-avoids both* $123$ *and* $321$*.*

Similarly, another result on key-avoiding ASMs can be easily deduced from pattern-avoidance in permutations.

**Proposition 2.3.** *For every set of patterns* $S$*, if* $231 \in S$*, then* $\mathrm{Av}^{\text{key}}(S) = \mathrm{pAv}(S)$*.*

**Proof:** If an ASM contains a $-1$ entry, then performing Step 2 of Definition 1.3 forces the appearance of a 231 pattern (with a "new" 1 playing the role of the "1" in 231 and the two 1s above and to the right of the $-1$, respectively, playing the roles of the "2" and the "3" in 231). Therefore, ASMs avoiding 231 are necessarily permutations, and the statement follows. □

| Excluded pattern | First terms of enumeration sequence | Known sequence |
|---|---|---|
| 1234 | 1, 2, 7, 41, 388, 5787 | Open |
| 1243 or 2134 | 1, 2, 7, 41, 393, 6085 | Open, see Prop. 2.7 for equivalence |
| 1324 | 1, 2, 7, 41, 394, 6161 | Open |
| 1342 or 2314 | 1, 2, 7, 40, 364, 5221 | Open, see Prop. 2.7 for equivalence |
| 1423 or 3124 | 1, 2, 7, 41, 393, 6084 | Open, see Prop. 2.7 for equivalence |
| 1432 or 3214 | 1, 2, 7, 41, 395, 6227 | Open, see Prop. 2.7 for equivalence |
| 2143 | 1, 2, 7, 41, 395, 6213 | Open |
| 2341 | 1, 2, 7, 37, 271, 2646 | Open |
| 2413 or 3142 | 1, 2, 7, 40, 362, 5156 | Open, see Prop. 2.7 for equivalence |
| 2431 or 3241 | 1, 2, 7, 39, 328, 4097 | Open, see Prop. 2.7 for equivalence |
| 3412 | 1, 2, 7, 38, 302, 3453 | Open |
| 3421 | 1, 2, 7, 39, 316, 3570 | Open |
| 4123 | 1, 2, 7, 41, 391, 5980 | Open |
| 4132 or 4213 | 1, 2, 7, 41, 392, 6036 | Open, see Prop. 2.7 for equivalence |
| 4231 | 1, 2, 7, 40, 354, 4720 | Open |
| 4312 | 1, 2, 7, 41, 384, 5594 | Open |
| 4321 | 1, 2, 7, 41, 386, 5671 | Open |

**Tab. 3:** Enumeration sequences of ASMs key-avoiding permutation patterns of size 4.

**Remark 2.4.** Note that $231$ key-avoiding ASMs are ASMs for *dominant* permutations. In the Schubert calculus literature, these are exactly the one term Schubert (and also Grothendieck) polynomials. So, from that perspective, it is immediate that there is a unique ASM with key $\pi$ if and only if $\pi$ avoids $231$ (and in this case that unique ASM is the permutation $\pi$).

**Corollary 2.5.** *For every $n$, $|\operatorname{Av}_n^{\text{key}}(231)| = \frac{1}{n+1}\binom{2n}{n}$.*

**Proof:** By Proposition 2.3, ASMs key-avoiding $231$ are exactly the permutations avoiding $231$, therefore counted by the Catalan numbers (see the foundational text [7, Section 2.2.1, Exercises 4 and 5]). □

**Corollary 2.6.** *For every $n$,*

$$|\operatorname{Av}_n^{\text{key}}(123, 231)| = \binom{n}{2} + 1, \text{ and}$$

$$|\operatorname{Av}_n^{\text{key}}(132, 231)| = |\operatorname{Av}_n^{\text{key}}(213, 231)| = |\operatorname{Av}_n^{\text{key}}(231, 312)| = |\operatorname{Av}_n^{\text{key}}(231, 321)| = 2^{n-1}.$$

**Proof:** As above, this follows combining Proposition 2.3 with the known enumeration of permutations avoiding these patterns, solved in [16]. □

For any permutation $\pi$, let us denote by $\pi^{\text{rev}}$ the reverse of the permutation $\pi$, which consists in reading its one-line notation from right to left.

**Proposition 2.7.** *For any permutation $\pi$, $\text{Av}^{\text{key}}(\pi)$ and $\text{Av}^{\text{key}}\left(\left((\pi^{rev})^{-1}\right)^{rev}\right)$ have the same enumeration sequence. This extends to sets of patterns $S$ and their image under $\pi \mapsto \left((\pi^{rev})^{-1}\right)^{rev}$.*

**Proof:** Observe that the function $\pi \mapsto \left((\pi^{\text{rev}})^{-1}\right)^{\text{rev}}$ corresponds to the reflection with respect to the northeast-southwest diagonal, with the convention of drawing permutation matrices as written just before Definition 1.3. This symmetry of the square is compatible with the computation of the key as we defined it (and it is actually the only one with this property). It therefore gives a size-preserving bijection between $\text{Av}^{\text{key}}(\pi)$ and $\text{Av}^{\text{key}}\left(\left((\pi^{\text{rev}})^{-1}\right)^{\text{rev}}\right)$, yielding the statement. □

Observing that the above symmetry maps 132 to 213, we have the following.

**Corollary 2.8.** *For every $n$, $|\text{Av}_n^{\text{key}}(132)| = |\text{Av}_n^{\text{key}}(213)|$.*

We also observe that whenever there is an enumeration coincidence in Table 2 or Table 3, it is explained by the symmetry identified by Proposition 2.7.

# 3 Monotone triangles and 312-key-avoidance

ASMs are in bijective correspondence with many other families of discrete objects (see e.g. [14]). One such family will play a particular role in our work: that of monotone triangles.

**Definition 3.1.** A *monotone triangle* of order $n$ is composed of $n$ rows, all left-aligned, the $i$-th row from the top containing $i$ elements, and such that

Condition (1): each row is strictly increasing from left to right;

Condition (2): each column is weakly increasing from bottom to top;

Condition (3): entries along each northwest-southeast diagonal are weakly increasing from top to bottom;

Condition (4): the bottom row is $1\,2\,\ldots\,n$.

**Definition 3.2.** As in [13, Section 5] (see also [2]), we associate a monotone triangle $t_A$ with any alternating sign matrix $A$ as follows:

Step (1): Create a 0-1 matrix $B$ by summing the rows of $A$ from the top to the current row.

Step (2): Build a monotone triangle by recording in its $i$-th row the indices of the columns where the 1s sit in row $i$ of matrix $B$.

**Example 3.3.** The process of obtaining a monotone triangle from an alternating sign matrix is shown below, on two examples.

$$A = \begin{pmatrix} 0 & 0 & 0 & 1 & 0 \\ 0 & 1 & 0 & -1 & 1 \\ 0 & 0 & 0 & 1 & 0 \\ 1 & -1 & 1 & 0 & 0 \\ 0 & 1 & 0 & 0 & 0 \end{pmatrix} \rightarrow B = \begin{pmatrix} 0 & 0 & 0 & 1 & 0 \\ 0 & 1 & 0 & 0 & 1 \\ 0 & 1 & 0 & 1 & 1 \\ 1 & 0 & 1 & 1 & 1 \\ 1 & 1 & 1 & 1 & 1 \end{pmatrix} \rightarrow t_A = \begin{array}{lllll} 4 \\ 2 & 5 \\ 2 & 4 & 5 \\ 1 & 3 & 4 & 5 \\ 1 & 2 & 3 & 4 & 5 \end{array}$$

$$A = \begin{pmatrix} 0 & 1 & 0 & 0 & 0 \\ 1 & -1 & 1 & 0 & 0 \\ 0 & 0 & 0 & 1 & 0 \\ 0 & 1 & 0 & -1 & 1 \\ 0 & 0 & 0 & 1 & 0 \end{pmatrix} \rightarrow B = \begin{pmatrix} 0 & 1 & 0 & 0 & 0 \\ 1 & 0 & 1 & 0 & 0 \\ 1 & 0 & 1 & 1 & 0 \\ 1 & 1 & 1 & 0 & 1 \\ 1 & 1 & 1 & 1 & 1 \end{pmatrix} \rightarrow t_A = \begin{array}{lllll} 2 \\ 1 & 3 \\ 1 & 3 & 4 \\ 1 & 2 & 3 & 5 \\ 1 & 2 & 3 & 4 & 5 \end{array}$$

The fact that this construction indeed produces a monotone triangle is well-known, but since it is central to our paper, we provide a short proof below.

**Proposition 3.4.** *The set of monotone triangles of order $n$ is in bijection with the set of $n \times n$ alternating sign matrices.*

**Proof:** Let $A$ be an $n \times n$ ASM and $t_A$ the output of the algorithm given above. We first show that the map is well-defined and the image is a monotone triangle. In Step (1) of the construction, we create a matrix $B$ of partial sums. It is important for Step (2) that all the entries of $B$ are in $\{0, 1\}$. By Definition 1.1, the nonzero entries of $A$ alternate between $1$ and $-1$ in each column, starting with a $1$. Thus, each column partial sum is either 0 or 1. Thus, Step (2) is well-defined, as is the entire construction.

We now show the image is a monotone triangle. So we need to show all the conditions of Definition 3.1 above are satisfied. Condition (1) is clear from the construction, since in Step (2), we recorded distinct numbers in increasing sequence. By the alternating condition (3. of Definition 1.1) and the fact that the rows and columns sum to 1 (2. of Definition 1.1), we have that the output of Step (2) is in the shape of a triangle. We also have Conditions (2) and (3) by the alternating condition (3. of Definition 1.1). Finally, Condition (4) is satisfied since each column of $A$ sums to 1.

We now note that each step of the map is bijective. Step (1) of the construction of $t_A$ is clearly invertible, since one may recover any matrix from its matrix of partial column sums. Step (2) is also invertible, since it merely records all 1 entries of the partial column sum matrix. □

Following [2], we say that a monotone triangle is *gapless* when in addition to the conditions of Definition 3.1, there are no "gaps" between column entries (in other words, the entries sitting in each column form an interval of integers). In Example 3.3, the first line shows a monotone triangle which is not gapless (3 is missing in the first column), and second line shows a gapless monotone triangle.

Ayyer, Cori, and Gouyou-Beauchamps [2] showed that the permutation matrices whose monotone triangles are gapless are exactly those that avoid 312.

**Proposition 3.5** ([2, Section 7])**.** *A permutation avoids the pattern* 312 *if and only if its associated monotone triangle is gapless.*

We now extend their result to show that 312-key-avoiding ASMs correspond exactly to gapless monotone triangles. We first need the following lemmas.

**Lemma 3.6.** *Suppose an ASM $A$ has a gapless monotone triangle and a $-1$ in position $(i,j)$ with the first $1$ to the west in its row in position $(i,j_0)$. Then there is no column index $j'$ with $j_0 < j' < j$ such that the first nonzero entry of $A$ south of row $i$ in column $j'$ (if it exists) equals $1$.*

**Proof:** Suppose $A$ is an $ASM$ with a $-1$ in position $(i,j)$ with the first 1 to its left in its row in position $(i,j_0)$. Suppose also there is a column index $j'$ with $j_0 < j' < j$ such that the first nonzero entry below row $i$ in column $j'$ equals 1, as shown below.

$$A = \begin{array}{c|ccccccccc} \vdots & \cdots & \vdots & \cdots & \vdots & \cdots & \vdots & \cdots \\ i & \cdots & 1 & 0\cdots0 & 0 & 0\cdots0 & -1 & \cdots \\ \vdots & \cdots & \vdots & \cdots & \vdots & \cdots & \vdots & \cdots \\ \vdots & \cdots & \cdots & \cdots & 1 & \cdots & \cdots & \cdots \\ \vdots & \cdots & \vdots & \cdots & \vdots & \cdots & \vdots & \cdots \\ \hline & \cdots & j_0 & \cdots & j' & \cdots & j & \cdots \end{array}$$

We compare rows $i-1$ and $i$ of the monotone triangle $t_A$. First, note that in $t_A$, $j$ is an entry in row $i-1$ but not in row $i$, while $j_0$ is not in row $i-1$ but does appear in row $i$. By our assumption on $j'$, observe that $j'$ does not appear in either of these two rows. Also, any other entries with values between $j_0$ and $j$ that appear in row $i-1$ also appear in row $i$; denote these entries as $c_1, c_2, \ldots c_\ell$ in increasing order. Finally, any nonzero entries in $A$ to the left of $(i,j_0)$ must come in pairs of a 1 followed by a $-1$, so rows $i-1$ and $i$ of $t_A$ have the same number of entries smaller than $j_0$. Thus, rows $i-1$ and $i$ of $t_A$ are as follows:

$$\begin{array}{ccccccc} \cdots & c_1 & c_2 & \cdots & c_\ell & j & \cdots \\ \cdots & j_0 & c_1 & \cdots & c_{\ell-1} & c_\ell & \cdots \end{array}$$

Since $j_0 < c_1 < c_2 < \ldots < c_\ell < j$, $j_0 < j' < j$, and $j'$ is distinct from the $c_1, \ldots, c_\ell$, there is some column in the above piece of the monotone triangle that has a gap, since $j'$ must be between two values in the sequence $j_0, c_1, \ldots, c_\ell, j$. Thus $t_M$ is not gapless, proving the lemma. □

Following Lemma 3.6, let us say that a $-1$ in any position $(i,j)$ of an ASM $A$ is *bad* when, letting $(i,j_0)$ be the position of the first 1 to its left in its row, there exists some column index $j'$ with $j_0 < j' < j$ such that the first nonzero entry of $A$ below row $i$ in column $j'$ exists and equals 1. This entry equal to 1 is called the *witness* that the $-1$ in position $(i,j)$ is bad. Non-bad $-1$ entries are referred to as *good.*

**Lemma 3.7.** *Suppose an ASM $A$ contains a bad $-1$ in position $(i,j)$. Let $A'$ be the ASM obtained after the removal of any other $-1$ in the SW key process (see Definition 1.3). Then the $-1$ entry in position $(i,j)$ is still bad in $A'$.*

**Proof:** As before, let us define $i_0$ such that the first 1 to the west of the $-1$ in position $(i,j)$ in its row is in position $(i,j_0)$. Because the $-1$ in position $(i,j)$ is bad, it has a corresponding witness, in column $j'$ for some $j_0 < j' < j$. Let us denote by $i' < i$ the row index of this witness. Without loss of generality (choosing the witness in the topmost possible row, and as far west as possible), we may assume that there is no 1 entry in the rectangle of positions $[i', i] \times [j_0+1, j']$ except for the witness. $A$ can therefore be represented as follows:

$$A = \begin{array}{c|ccccccc}
\vdots & \cdots & \vdots & \cdots & \vdots & \cdots & \vdots & \cdots \\
i_0 & \cdots & \cdots & \cdots & -1 & \cdots & \cdots & \cdots \\
\vdots & \cdots & \vdots & \cdots & \vdots & \cdots & \vdots & \cdots \\
i & \cdots & 1 & 0\cdots 0 & 0 & 0\cdots 0 & -1 & \cdots \\
\vdots & \cdots & \vdots & \cdots & \vdots & \cdots & \vdots & \cdots \\
i' & \cdots & 0 & 0\cdots 0 & 1 & \cdots & \cdots & \cdots \\
\vdots & \cdots & \vdots & \cdots & \vdots & \cdots & \vdots & \cdots \\
\hline
 & \cdots & j_0 & \cdots & j' & \cdots & j & \cdots
\end{array}.$$

During the removal of any $-1$, any $1$ entry that is affected by this removal (except the one in the column of the removed $-1$ and below it) yields a $1$ entry in the ASM after the removal, which is in the same column but lower. So, the only possible removal of a $-1$ in $A$ that could turn the bad $-1$ in position $(i, j)$ into a good $-1$ would be the removal of a $-1$ in position $(i_0, j')$ for some $i_0 > i$ (and such that all entries in positions $(k, j')$ for $i' < k < i_0$ are $0$). In particular, this requires that this $-1$ in position $(i_0, j')$ is removable.

Consequently, the rectangle of positions $[i', i] \times [j_0 + 1, j']$ may not contain any $-1$ entries. And since we have already seen that it does not contain $1$ entries except for the witness, all entries there are $0$.

During the removal of the identified $-1$ in position $(i_0, j')$, two cases may occur: denoting as $j'_0$ the column of the closest $1$ to its left in its row, we may have $j'_0 > j_0$ or $j'_0 < j_0$ (equality is impossible, since the $-1$ in position $(i_0, j')$ is removable).

- If $j'_0 > j_0$, then, after the removal of the $-1$ in position $(i_0, j')$, there will be a $1$ entry in row $i'$, in some column $j'_0$ and $j' - 1$, which serves as a witness that the $-1$ in position $(i, j)$ is still bad after the removal. This is illustrated in Figure 2.

- The situation when $j'_0 < j_0$ is similar, the only difference being that the closest $1$ to the west of the $-1$ in position $(i, j)$ may have shifted further to the west. Whether or not this is the case, there is an entry $1$ that serves as a witness that the $-1$ in position $(i, j)$ is still bad after the removal. This is illustrated in Figure 3.

□

We now prove the main theorem of this section.

**Theorem 3.8.** *The monotone triangle construction $A \mapsto t_A$ described above puts in bijection* 312*-key-avoiding $n \times n$ ASMs and gapless monotone triangles with bottom row* $1\ 2\ \ldots\ n$.

**Proof:** We wish to show that the key of an ASM $A$ avoids 312 if and only if the monotone triangle of $A$ is gapless. We shall do so by induction on the number of negative ones in an ASM.

First observe that our statement is true when restricted to ASMs with no negative ones. Indeed, such an ASM is its own key, and Proposition 3.5 above (from [2]) proves that the permutation matrices whose monotone triangles are gapless are exactly the 312-avoiding permutations.

$$
\begin{array}{c|ccccccc}
i_0 & \cdots & \cdots & 1 & 0\cdot\cdot 0 & -1 & \cdots & \cdots \\
\vdots & \vdots & \cdots & \vdots & \cdots & \vdots & \cdots & \vdots \\
i & 1 & 0\cdot\cdot 0 & 0 & 0\cdot\cdot 0 & 0 & 0\cdot\cdot 0 & -1 \\
\vdots & \vdots & \cdots & \vdots & \cdots & \vdots & \cdots & \vdots \\
i' & 0 & 0\cdot\cdot 0 & \vdots & 0\cdot\cdot 0 & \boxed{1} & \cdots & \cdots \\
\vdots & \vdots & \cdots & \vdots & \cdots & \vdots & \cdots & \vdots \\
\hline
 & j_0 & \cdots & j_0' & \cdots & j' & \cdots & j
\end{array}
\rightarrow
\begin{array}{c|ccccccc}
i_0 & \cdots & \cdots & 0 & 0\cdot\cdot 0 & -1 & \cdots & \cdots \\
\vdots & \vdots & \cdots & \vdots & \cdots & \vdots & \cdots & \vdots \\
i & 1 & 0\cdot\cdot 0 & 0 & 0\cdot\cdot 0 & 0 & 0\cdot\cdot 0 & -1 \\
\vdots & \vdots & \cdots & \vdots & \cdots & \vdots & \cdots & \vdots \\
i' & 0 & 0\cdot\cdot 0 & \vdots & \cdot\boxed{1}\cdot & 0 & \cdots & \cdots \\
\vdots & \vdots & \cdots & \vdots & \cdots & \vdots & \cdots & \vdots \\
\hline
 & j_0 & \cdots & j_0' & \cdots & j' & \cdots & j
\end{array}
$$

**Fig. 2:** How a new witness is formed if $j_0 < j_0' < j$.

$$
\begin{array}{c|ccccccc}
i_0 & 1 & 0\cdot\cdot 0 & 0 & 0\cdot\cdot 0 & -1 & \cdots & \cdots \\
\vdots & \vdots & \cdots & \vdots & \cdots & \vdots & \cdots & \vdots \\
i & \cdots & \cdots & 1 & 0\cdot\cdot 0 & 0 & 0\cdot\cdot 0 & -1 \\
\vdots & \vdots & \cdots & \vdots & \cdots & \vdots & \cdots & \vdots \\
i' & 0 & 0\cdot\cdot 0 & \vdots & 0\cdot\cdot 0 & \boxed{1} & \cdots & \cdots \\
\vdots & \vdots & \cdots & \vdots & \cdots & \vdots & \cdots & \vdots \\
\hline
 & j_0' & \cdots & j_0 & \cdots & j' & \cdots & j
\end{array}
\rightarrow
\begin{array}{c|ccccccc}
i_0 & 0 & 0\cdot\cdot 0 & 0 & 0\cdot\cdot 0 & -1 & \cdots & \cdots \\
\vdots & \vdots & \cdots & \vdots & \cdots & \vdots & \cdots & \vdots \\
i & \cdot\cdot & \cdot 1\cdot & 0 & 0\cdot\cdot 0 & 0 & 0\cdot\cdot 0 & -1 \\
\vdots & \vdots & \cdots & \vdots & \cdots & \vdots & \cdots & \vdots \\
i' & 0 & 0\cdot\cdot 0 & 0 & \cdot\boxed{1}\cdot & 0 & \cdots & \cdots \\
\vdots & \vdots & \cdots & \vdots & \cdots & \vdots & \cdots & \vdots \\
\hline
 & j_0' & \cdots & \cdots & \cdots & j' & \cdots & j
\end{array}
$$

**Fig. 3:** How a new witness is formed if $j_0' < j_0$.

Suppose now that our statement is true for alternating sign matrices with at most $k-1$ negative ones, and proceed by proving separately each of the two implications in our claim for $k \geq 1$ negative ones.

So, suppose first that $A$ is an ASM with $k$ negative ones whose monotone triangle $t_A = M$ is gapless. We need to show that removing a negative one from $A$ as prescribed by the process of computing the key does not introduce a gap in the resulting monotone triangle $M'$. Pick a $-1$ that has no other $-1$ southwest of it. Let $(i, j)$ be the position of this $-1$, $(i, j_0)$ be the position of the 1 to its left in its row, and $(i_0, j)$ be the position of the 1 below it in its column. Let $A'$ denote the ASM produced by removing this $-1$ by the SW key process, and $M'$ its corresponding monotone triangle.

By Lemma 3.6, there is no column index $j'$ with $j_0 < j' < j$ such that the first nonzero entry of $A$ below row $i$ in column $j'$ equals 1. So the removal of this $-1$ by the SW-key process is a simple removal:

$$\begin{matrix} 1 & \dots & -1 \\ \vdots & \dots & \vdots \\ 0 & \dots & 1 \end{matrix} \quad \longrightarrow \quad \begin{matrix} 0 & \dots & 0 \\ \vdots & \dots & \vdots \\ 1 & \dots & 0 \end{matrix}$$

where the regions indicated by ellipses contain only zeros.

Since the change in the matrix is isolated, we can directly compare the portions of the monotone triangle from row $i$ to row $i_0$ of $M$ (on the left below) with the same rows of $M'$ (on the right below):

$$\begin{matrix} \cdots j_0 & j_0+1 & \dots & j-2 & j-1 & \ell & \cdots \\ \vdots & \vdots & \vdots & \vdots & \vdots & \vdots & \\ \cdots j_0 & j_0+1 & \dots & j-2 & j-1 & \ell & \cdots \\ \cdots j_0 & j_0+1 & \dots & j-2 & j-1 & j & \cdots \end{matrix} \quad \longrightarrow \quad \begin{matrix} \cdots & j_0+1 & j_0+2 & \dots & j-1 & j & \ell\cdots \\ & \vdots & \vdots & \vdots & \vdots & \vdots & \vdots \\ \cdots & j_0+1 & j_0+2 & \dots & j-1 & j & \ell\cdots \\ \cdots & j_0 & j_0+1 & \dots & j-2 & j-1 & j\cdots \end{matrix}$$

for some $\ell > j$. Since $M$ is gapless, $\ell = j+1$. The above comparison of $M$ and $M'$ shows that the same region of $M'$ is gapless. $M$ and $M'$ agree outside this region, so $M'$ is gapless. By the induction hypothesis, since $M'$ is gapless and $A'$ has $k-1$ negative ones, the key of $A'$ avoids 312. But the key of $A$ equals the key of $A'$, so $A$ is 312-key-avoiding, proving the first part of our claim.

Now, let us suppose that $A$ is an ASM with $k$ negative ones whose key avoids 312. We wish to show its monotone triangle $M$ is gapless. Since $A$ key-avoids 312, then removing any $-1$ according to the SW-key process produces an ASM $A'$ with $k-1$ negative ones that has the same key. So by our induction hypothesis, since the key of $A'$ avoids 312, the monotone triangle $M'$ of $A'$ is gapless.

Since $M'$ is gapless, Lemma 3.6 ensures that $A'$ has no bad $-1$s. Let us first prove by contradiction that $A$ has no bad $-1$. By Lemma 3.7, if $A$ has a bad $-1$, it must have been the one we removed to obtain $A'$. Thus, the bad $-1$ must be removable, so there is no other $-1$ southwest of it. It is also the only removable $-1$ of $A$, since otherwise, we could remove the other $-1$ first and obtain some $A'$ which would still have a bad $-1$ (again by Lemma 3.7).

So let $(i, j)$ be the position of this bad $-1$ in $A$, assumed to exist. Then let the closest 1 to its left in its row be in position $(i, j_0)$ and the closest 1 below it in its column be in position $(i_0, j)$. Then let $(i', j')$, with $j_0 < j' < j$, be such that column $j'$ has its first nonzero entry below row $i$ equals 1, and chosen such that $i'$ is minimal among all possible values of $i'$. We compare rows $i_0 - 1$ and $i_0$ of $M$ and $M'$ in two cases (the third case $i' = i_0$ being impossible since the considered $-1$ is removable).

Case $i' < i_0$: In this case, the $-1$ removal is not simple. We compare rows $i-1$ and $i$ of $M$ with rows $i'-1$ and $i'$ of $M'$. As in the proof of Lemma 3.6, let $c_1, \dots, c_\ell$ denote the column indices strictly between $j_0$ and $j$ which have a column partial sum of 1 in row $i$. Note, there is at least one gap in this list at $j'$. So rows $i-1$ and $i$ of $M$ are as follows, and there is some gap:

$$\begin{array}{cccccc} \cdots & c_1 & c_2 & \ldots & c_\ell & j \\ \cdots & j_0 & c_1 & \ldots & c_{\ell-1} & c_\ell \end{array}$$

We now look at rows $i'-1$ and $i'$ of $M'$ and note they are exactly the same as the above:

$$\begin{array}{cccccc} \cdots & c_1 & c_2 & \ldots & c_\ell & j \\ \cdots & j_0 & c_1 & \ldots & c_{\ell-1} & c_\ell \end{array}$$

So $M'$ must also have a gap; this is a contradiction.

Case $i' > i_0$: In this case, we see that the rows $i_0-1$ and $i_0$ of $M'$ are the same as shown above for rows $i-1$ and $i$ of $M$, deriving a similar contradiction. Thus, there is no bad $-1$ in $M$.

Finally, we need to consider the case that $A$ has no bad $-1$, and $A'$ is obtained by removing a good (and removable) $-1$ from $A$. As before, denote by $(i,j)$ the position of the removed $-1$, and let the closest 1 to its left in its row be in position $(i,j_0)$ and the closest 1 below it in its column be in position $(i_0,j)$. Our current assumption means there is no $j'$ between $j_0$ and $j$ such that the first nonzero entry in that column of $A$ below row $i$ is a 1. So all numbers strictly between $j_0$ and $j$ appear in rows $i-1$ and $i$ of both $M$ and $M'$. We compare rows $i-1$ through $i_0$ of $M$ and $M'$.

$$\begin{array}{cccccc} \cdots j_0+1 & j_0+2 & \cdots & j-1 & j & \ell\cdots \\ \cdots j_0 & j_0+1 & \cdots & j-2 & j-1 & \ell\cdots \\ \vdots & \vdots & \vdots & \vdots & \vdots & \vdots \\ \cdots j_0 & j_0+1 & \cdots & j-2 & j-1 & \ell\cdots \\ \cdots j_0 & j_0+1 & \cdots & j-2 & j-1 & j\cdots \end{array} \longrightarrow \begin{array}{cccccc} \cdots j_0+1 & j_0+2 & \cdots & j-1 & j & \ell\cdots \\ \cdots j_0+1 & j_0+2 & \ldots & j-1 & j & \ell\cdots \\ \vdots & \vdots & \vdots & \vdots & \vdots & \vdots \\ \cdots j_0+1 & j_0+2 & \ldots & j-1 & j & \ell\cdots \\ \cdots j_0 & j_0+1 & \ldots & j-2 & j-1 & j\cdots \end{array}$$

Now since $A$ is 312-key-avoiding, $A'$ is also, since they have the same key. By the induction hypothesis, since $A'$ is 312-key-avoiding and $A'$ has $k-1$ negative ones, $M'$ is gapless. Thus it must be that $\ell = j+1$. This implies $M$ also has no gaps on the region shown above. Since $M$ is equal to $M'$ outside this region, we have that $M$ is gapless, concluding the proof. □

**Remark 3.9.** In [2, Sec. 9], a notion of 312-pattern-avoidance is given on *gog words*, objects the authors define that are in bijection with monotone triangles. They prove this is equivalent to the gapless property of a monotone triangle. Thus, it follows from Theorem 3.8 that their notion of 312-pattern-avoidance corresponds to SW 312-key-avoidance. It is noted in [2, Sec. 10] that it was unclear how to extend their notion of 312-avoidance on gog words to other patterns. Thus, key-avoidance may be seen as a natural generalization of their notion of 312-avoidance that works for all patterns.

# 4 Catalan numbers and key-avoidance of 312 and 321

In this section, we prove that $\mathrm{Av}^{\mathrm{key}}(312, 321)$ is enumerated by the Catalan numbers, see Theorem 4.8. The proof uses the following series of lemmas. We begin by noting that the 1 entries affected by a given removal in the SW key process must form an increasing pattern (in the sense of classical patterns).

**Lemma 4.1.** *Consider an ASM $A$ and a removable $-1$ in $A$. Denote by $A'$ the ASM obtained from $A$ after the removal of this $-1$ by the SW key process. Consider also the $m+1$ corresponding $1$s in $A$ that are moved during this removal. These $1$s form an occurrence of the classical increasing pattern of size $m+1$ in $A$, which is mapped by the considered removal to an occurrence of the classical increasing pattern of size $m$ in $A'$, where the largest (column) value of the original increasing sequence is removed and all other (column) values are moved to larger row indices.*

**Proof:** Any new 1 appearing in position $(i, j)$ in $A'$ must result from a 1 in position $(i-a, j)$ and another 1 in position $(i, j+b)$ in $A$, where $a, b > 0$. Denote by $(i_0, j_0), (i_1, j_1), \dots, (i_m, j_m)$ the positions of the 1s in $A$ affected by the considered removal, in order from top to bottom in $A$. That is, $i_0 < i_1 < \cdots < i_m$ and $j_0 < j_1 < \cdots < j_m$. Then this removal switches the 1s in the original sequence to 0s and creates a new sequence of 1s in positions $(i_1, j_0), (i_2, j_1), \dots, (i_m, j_{m-1})$. Thus the largest (column) value of the increasing sequence $j_m$ is removed and each smaller value moves forward to the original (row) index of the next term. □

**Example 4.2.** See Figure 4, with the removal of the $-1$ at $(2, 4)$ by the SW key process. The affected 1s are in positions $\{(2,1),(3,3),(4,4)\}$ and form an increasing sequence. After the removal, new 1s are in positions $\{(3,1),(4,3)\}$, also forming an increasing sequence, of size one less, and in the 'shifted' positions prescribed in the proof of Lemma 4.1.

$$\begin{pmatrix} 0 & 0 & 0 & 1 & 0 \\ \boxed{1} & 0 & 0 & -1 & 1 \\ 0 & 0 & \boxed{1} & 0 & 0 \\ 0 & 0 & 0 & \boxed{1} & 0 \\ 0 & 1 & 0 & 0 & 0 \end{pmatrix} \to \begin{pmatrix} 0 & 0 & 0 & 1 & 0 \\ 0 & 0 & 0 & 0 & 1 \\ \boxed{1} & 0 & 0 & 0 & 0 \\ 0 & 0 & \boxed{1} & 0 & 0 \\ 0 & 1 & 0 & 0 & 0 \end{pmatrix}$$

**Fig. 4:** The SW key process taking the classical increasing pattern 123 occurring at entries $\{(2,1),(3,3),(4,4)\}$ to the classical increasing pattern 12 at entries $\{(3,1),(4,3)\}$.

As a consequence of Lemma 4.1, we can then show that no (classical) 321 pattern can be formed during the SW key process without there already having been a (classical) 321 pattern.

**Lemma 4.3.** *If $A$ is an ASM that classically avoids* 321*, then the SW key of $A$ also avoids* 321*.*

**Proof:** Let $A$ be an ASM that classically avoids 321. If by way of contradiction, the key of $A$ were to classically contain 321, then at some removal of a $-1$ by the SW key process, a 321 pattern must be formed. Let us denote by $B$ (resp. $C$) the ASM just before (resp. just after) the removal of the $-1$ which creates the 321 pattern.

By Lemma 4.1, we know that the removal of a $-1$ moves an increasing sequence into an increasing sequence. Therefore, exactly one entry of the newly formed 321 pattern must have been affected by this removal. Thus, the ASM $B$ must satisfy the following.

1. It must have a decreasing pattern of size two, say in entries $(x_1, d), (x_3, b)$ with $x_1 < x_3$ and $d > b$, that is unaffected by the current removal of a $-1$.

2. The removal of this $-1$ must change the (row) index of an entry in the affected increasing sequence such that it becomes part of a decreasing pattern with $(x_1, d), (x_3, b)$.

First note that the 1 entries in positions $(x_1, d)$ and $(x_3, b)$ cannot act as 2 and 1 in a 321 pattern in $C$. Indeed, if they did, there would be a 1 entry to the NE of $(x_1, d)$ in $C$, which belongs to the increasing sequence affected by the removal of the $-1$; but this 1 entry of $C$ has a corresponding 1 entry in $B$ to its NE, contradicting that $B$ avoids 321.

Hence the removal of the $-1$ must satisfy one or the other of the following.

Case 1: This $-1$ removal moves some 1 entry in position $(x_0, c)$ in $B$ to position $(x_2, c)$ in $C$, for some $x_0 < x_1 < x_2 < x_3$ and $b < c < d$.

Case 2: This $-1$ removal moves some 1 entry in position $(x_\alpha, a)$ in $B$ to position $(x_4, a)$ in $C$, for some $x_\alpha \leq x_3 < x_4$ and $a < b < d$.

In Case 1 (illustrated on Fig. 5, left), the 1 in position $(x_0, c)$ in $B$ would be moving due to the removal of the $-1$ to the (row) index of a 1 in position $(x_2, e)$ in $B$, as described in the proof of Lemma 4.1. This larger (column) value $e$ thus must be such that $e > d$ for there not to have already been a 321 pattern formed by 1 entries at positions $(x_1, d)$, $(x_2, e)$ and $(x_3, b)$. But then the removed $-1$ must be to the northeast of $(x_1, d)$. However, this cannot happen as the 1 in position $(x_1, d)$ or another 1 entry to the northeast of $(x_1, d)$ would have also been part of the 1 entries affected by the considered $-1$ removal, and thus would have blocked this move from happening.

In Case 2 (illustrated on Fig. 5, right), we have a parallel situation with the 1 entry of $B$ in position $(x_\alpha, a)$ where $x_\alpha \leq x_3$ takes the (row) index from a 1 in position $(x_4, f)$ in $B$, with $x_4 > x_3$. Note that $f > b$ for there not to have already been a 321 pattern in $B$ formed by 1 entries at positions $(x_1, d)$, $(x_3, b)$ and $(x_4, f)$. Now one can see that the removed $-1$ is to the northeast of the 1 in position $(x_3, b)$. Thus this 1 in position $(x_3, b)$ or another 1 to the northeast of $(x_3, b)$ must be part of the 1 entries affected by the considered $-1$ removal, hence blocking this move from happening. □

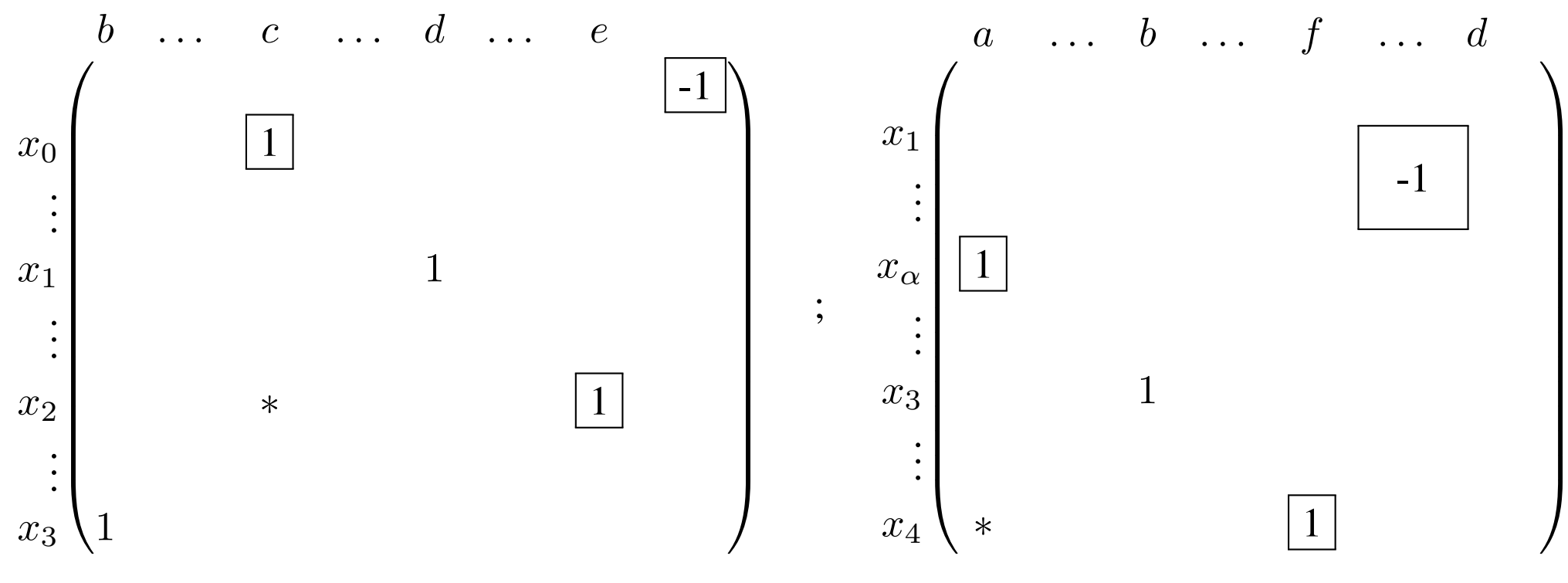

**Fig. 5:** Illustrating Cases 1 (left) and 2 (right) of the proof of Lemma 4.3. The boxed 1 entries are those affected by the removal of the considered $-1$, which is also boxed. The position marked with $*$ indicates the position (in $C$) to which the 1 entry in position $(x_0, c)$ (resp. $(x_\alpha, a)$) of $B$ is moved.

In fact, the proof above gives us the following stronger result.

**Corollary 4.4.** *If $A$ is an ASM that classically avoids* 321*, then any ASM obtained at some stage during the computation of the SW key of $A$ also classically avoids* 321*.*

For an ASM whose SW key avoids 312, the converse of Lemma 4.3 is true as well.

**Lemma 4.5.** *Let $A$ be an ASM whose SW key avoids* 312*. If $A$ classically contains* 321*, then the SW key of $A$ also contains* 321*.*

**Proof:** Recall that the SW key of an ASM is not dependent on the order in which we remove removable $-1$s. Without loss of generality, in this proof, we therefore consider the deterministic process that consists in removing the south-most removable $-1$ whenever we have a choice of removable $-1$.

Let us prove the following claim, which immediately implies the announced statement: the removal of the south-most removable $-1$ from an ASM classically containing 321 results in an ASM that still classically contains 321 or contains a 312 pattern that cannot be removed later by the SW key process.

Throughout the proof, for an ASM $M$ where an occurrence of a 321 or 312 pattern is identified, when referring to the smallest (resp. largest, resp. middle) entry of this pattern in $M$, we mean the entry 1 in $M$ that corresponds to the element 1 (resp. 3, resp. 2) of the pattern.

Let us first prove that the removal of a removable $-1$ from an ASM classically containing 321 results in an ASM that still classically contains 321 or contains a 312 pattern. To this effect, let us consider two ASMs $B$ and $C$ such that $B$ classically contains 321 and $C$ is obtained from $B$ by the removal of a removable $-1$.

Let $(x_1, y_3), (x_2, y_2), (x_3, y_1)$ be a 321 pattern of 1s in $B$, i.e. $x_1 < x_2 < x_3$ and $y_3 > y_2 > y_1$. Since Lemma 4.1 ensures that the 1 entries of $B$ affected by the considered removal form an increasing sequence, it follows that exactly one of these three entries of $B$ is affected by this removal.

1. Assume that the smallest entry $(x_3, y_1)$ of the considered 321 pattern in $B$ is the one affected by the considered removal. The SW key process forces the removal to result in a new 1 to the west (smaller value) and/or to the south (later position) of $(x_3, y_1)$. Thus there exists a 1 entry in $C$ in at least one of the positions $(x_3, y_0), (x_4, y_1)$ where $y_0 < y_1, x_3 < x_4$, and it can play the role of the "1" in a 321 pattern where the "2" and "3" remain the 1 entries in $(x_2, y_2)$ and $(x_1, y_3)$ (unchanged from $B$ to $C$). An example of this is shown in Figure 6. This shows that $C$ contains 321.

$$\underbrace{\begin{pmatrix} 0 & 0 & 0 & 0 & \boxed{1} \\ 0 & 0 & \boxed{1} & 0 & 0 \\ 1 & 0 & -1 & 1 & 0 \\ 0 & \boxed{1} & 0 & 0 & 0 \\ 0 & 0 & 1 & 0 & 0 \end{pmatrix}}_{B} \rightarrow \underbrace{\begin{pmatrix} 0 & 0 & 0 & 0 & \boxed{1} \\ 0 & 0 & \boxed{1} & 0 & 0 \\ 0 & 0 & 0 & 1 & 0 \\ \boxed{1} & 0 & 0 & 0 & 0 \\ 0 & \boxed{1} & 0 & 0 & 0 \end{pmatrix}}_{C}$$

**Fig. 6:** When the smallest entry of a classical 321 pattern is affected by the removal of a $-1$ in the SW key process, there is at least one new 1 to take its place.

2. Assume next that the largest entry $(x_1, y_3)$ of the considered 321 pattern in $B$ is the one affected by the considered removal. In particular, the removed $-1$ is weakly to the northeast of $(x_1, y_3)$, and has $1$ entries above it in its column and to its right in its row. Either of these $1$ entries can play the role of the "3" in a $321$ pattern in $C$, where the "1" and "2" remain the $1$ entries in $(x_3, y_1)$ and $(x_2, y_2)$ (unchanged from $B$ to $C$), as shown in Figure 7. Again, this shows that $C$ contains 321.

$$\underbrace{\begin{pmatrix} 0 & 0 & 0 & 1 & 0 \\ 0 & 0 & \boxed{1} & -1 & 1 \\ 0 & \boxed{1} & 0 & 0 & 0 \\ 0 & 0 & 0 & 1 & 0 \\ \boxed{1} & 0 & 0 & 0 & 0 \end{pmatrix}}_{B} \rightarrow \underbrace{\begin{pmatrix} 0 & 0 & 0 & \boxed{1} & 0 \\ 0 & 0 & 0 & 0 & \boxed{1} \\ 0 & \boxed{1} & 0 & 0 & 0 \\ 0 & 0 & 1 & 0 & 0 \\ \boxed{1} & 0 & 0 & 0 & 0 \end{pmatrix}}_{C}$$

**Fig. 7:** When the largest entry of a classical $321$ pattern is affected by the removal of a $-1$ in the SW key process, either one of the $1$ to the north and the $1$ to the east of the removed $-1$ can be used to form another $321$ pattern with the original $21$ part of the pattern.

3. The only case left to consider is when the $1$ entry of $B$ that is affected by the considered removal is $(x_2, y_2)$. Let us show that if $C$ classically avoids 321, then it must contain 312. This follows easily, since for $C$ to classically avoid 321, the new 1(s) in $C$ are all to the west or south of $(x_3, y_1)$ (but not both), implying that at least one 312 pattern is formed in $C$. Note that the $1$ in position $(x_3, y_1)$ could serve as either the $1$ (as shown in Figure 8, first row) or the $2$ (as shown in Figure 8, second row) of the 312 pattern depending on the placement of the removed $-1$ (and the new $1$ in $C$ would take on the other role).

$$\underbrace{\begin{pmatrix} 0 & 0 & \boxed{1} & 0 & 0 \\ 0 & \boxed{1} & -1 & 1 & 0 \\ 0 & 0 & 0 & 0 & 1 \\ \boxed{1} & 0 & 0 & 0 & 0 \\ 0 & 0 & 1 & 0 & 0 \end{pmatrix}}_{B} \rightarrow \underbrace{\begin{pmatrix} 0 & 0 & \boxed{1} & 0 & 0 \\ 0 & 0 & 0 & 1 & 0 \\ 0 & 0 & 0 & 0 & 1 \\ \boxed{1} & 0 & 0 & 0 & 0 \\ 0 & \boxed{1} & 0 & 0 & 0 \end{pmatrix}}_{C}$$

or

$$\underbrace{\begin{pmatrix} 0 & 0 & 1 & 0 & 0 \\ 1 & 0 & -1 & \boxed{1} & 0 \\ 0 & 0 & \boxed{1} & 0 & 0 \\ 0 & \boxed{1} & 0 & 0 & 0 \\ 0 & 0 & 0 & 0 & 1 \end{pmatrix}}_{B} \rightarrow \underbrace{\begin{pmatrix} 0 & 0 & 1 & 0 & 0 \\ 0 & 0 & 0 & \boxed{1} & 0 \\ \boxed{1} & 0 & 0 & 0 & 0 \\ 0 & \boxed{1} & 0 & 0 & 0 \\ 0 & 0 & 0 & 0 & 1 \end{pmatrix}}_{C}$$

**Fig. 8:** When the middle entry of a classical $321$ pattern is affected by the removal of a $-1$ in the SW key process and the $321$ is destroyed, so a classical $312$ must be formed.

This concludes the first part of our proof. We now turn to proving that, if a 312 has been created as in case 3. above during the SW key process, then the SW key process will always result in an ASM classically containing 312.

Let us therefore consider an ASM $C$ classically containing a 312 pattern that has been formed from a 321 pattern due to the removal(s) of one (or several) $-1$, and let $D$ be the ASM obtained after removing the south-most removable $-1$ from $C$. We want to show that $D$ also contains a classical 312 pattern. More precisely, we show that $D$ contains such a pattern with the two smaller entries of the pattern in the same rows as in $C$. (This allows the arguments below to also apply when $C$ has been obtained several removals after the first appearance of a 312 pattern.)

Because the 3 forms a decreasing sequence with each of the $1, 2$ in a 312 pattern, by Lemma 4.1, no single removal of a $-1$ will affect both the largest entry and also one of the smaller two entries of a 312 pattern. Next observe that if the removal affects the 1 entry of $C$ corresponding to the largest entry of a 312 pattern in $C$, then $C$ must have a $-1$ weakly to the northeast of this 1 entry. Thus, there exists a 1 to the north and another to the east of that $-1$, either of which can replace the original 3 of the 312 pattern in $C$, yielding a 312 pattern in $C$, using the same two smaller entries as the 312 pattern we started from in $C$, and that stays the same in $D$. This is shown in Figure 9.

$$\underbrace{\begin{pmatrix} 0 & 0 & 0 & 1 & 0 \\ 0 & 0 & \boxed{1} & -1 & 1 \\ \boxed{1} & 0 & 0 & 0 & 0 \\ 0 & \boxed{1} & 0 & 0 & 0 \\ 0 & 0 & 0 & 1 & 0 \end{pmatrix}}_{C} \rightarrow \underbrace{\begin{pmatrix} 0 & 0 & 0 & \boxed{1} & 0 \\ 0 & 0 & 0 & 0 & \boxed{1} \\ \boxed{1} & 0 & 0 & 0 & 0 \\ 0 & \boxed{1} & 0 & 0 & 0 \\ 0 & 0 & 1 & 0 & 0 \end{pmatrix}}_{D}$$

**Fig. 9:** When the largest entry of a classical 312 pattern is affected by the removal of a $-1$ in the SW key process in $C$, we still find a 312 pattern in $D$, with the two smaller entries unchanged.

Thus consider the removal of a $-1$ that affects at least one of the smaller entries of the pattern 312 in $C$. Denote by $(a_1, b_1), (a_2, b_2)$ with $a_1 < a_2$ and $b_1 < b_2$ the positions of these entries in $C$. Let us also denote by $(a_0, b_3)$ the position of the largest entry of the considered 312 pattern in $C$. In particular, $a_0 < a_1$ and $b_2 < b_3$. The $-1$ in $C$ that we remove to get $D$ is weakly north of the previously removed $-1$(s) by our stated convention noted at the beginning of the proof. Thus this $-1$ is strictly north of both 1s representing the $1, 2$ of the 312 pattern in $C$. (Fig. 8 can help visualize the above claim.) Because of this, whether the removal affects one or both entries, there must be at least one more 1 in the Ferrer's board affected by the removal of the $-1$: the one located west of the $-1$ in its row, and which will have the smallest row index of all affected 1s.

1. If the 1 in both positions $(a_1, b_1), (a_2, b_2)$ are affected, then since neither has the smaller (row) index of the set of affected 1s, by (the proof of) Lemma 4.1, we must have new 1s in positions $(a_1, c_1), (a_2, c_2)$ in $D$ where $c_1 < b_1$, $c_2 < b_2$, and $c_1 < c_2$. Hence the entries in positions $(a_0, b_3), (a_1, c_1), (a_2, c_2)$ form a 312 pattern in $D$ as shown in Figure 10.

2. Suppose now only the 1 in position $(a_2, b_2)$ in $C$ is affected. The $-1$ to be removed is east of $(a_1, b_1)$ because $b_1 < b_2$. Thus there must be another 1 entry in position $(c, d)$ east of $(a_1, b_1)$ (i.e. $b_1 < d$), which is the next 1 northwest of $(a_2, b_2)$ (i.e. $c < a_2$ and $d < b_2$) in the affected Ferrer's

$$\underbrace{\begin{pmatrix} 0 & 0 & 0 & \boxed{1} & 0 \\ 1 & 0 & 0 & -1 & 1 \\ 0 & \boxed{1} & 0 & 0 & 0 \\ 0 & 0 & \boxed{1} & 0 & 0 \\ 0 & 0 & 0 & 1 & 0 \end{pmatrix}}_{C} \to \underbrace{\begin{pmatrix} 0 & 0 & 0 & \boxed{1} & 0 \\ 0 & 0 & 0 & 0 & 1 \\ \boxed{1} & 0 & 0 & 0 & 0 \\ 0 & \boxed{1} & 0 & 0 & 0 \\ 0 & 0 & 1 & 0 & 0 \end{pmatrix}}_{D}$$

**Fig. 10:** When the two smaller entries of a classical 312 pattern (formed by the removal of a 321) are both affected by the removal of a $-1$ in the SW key process, a 312 is formed in the same rows as the original 312 pattern.

board. Thus a new $1$ is formed in position $(a_2, d)$ in $D$ and so $(a_0, b_3), (a_1, b_1), (a_2, d)$ form a 312 pattern in $D$ with the two smaller entries in rows $a_1$ and $a_2$, as shown in Figure 11.

$$\underbrace{\begin{pmatrix} 0 & 0 & 0 & \boxed{1} & 0 \\ 0 & 1 & 0 & -1 & 1 \\ \boxed{1} & 0 & 0 & 0 & 0 \\ 0 & 0 & \boxed{1} & 0 & 0 \\ 0 & 0 & 0 & 1 & 0 \end{pmatrix}}_{C} \to \underbrace{\begin{pmatrix} 0 & 0 & 0 & \boxed{1} & 0 \\ 0 & 0 & 0 & 0 & 1 \\ \boxed{1} & 0 & 0 & 0 & 0 \\ 0 & \boxed{1} & 0 & 0 & 0 \\ 0 & 0 & 1 & 0 & 0 \end{pmatrix}}_{D}$$

**Fig. 11:** When only the middle entry of a classical 312 pattern (formed by the removal of a 321) is affected by the removal of a $-1$ in the SW key process, a 312 is formed in the same rows as the original pattern.

3. Finally, if only the $1$ in position $(a_1, b_1)$ is affected, again by (the proof of) Lemma 4.1, a new $1$ is formed in position $(a_1, b_0)$ in $D$ for some $b_0 < b_1$. Thus $(a_0, b_3), (a_1, b_0), (a_2, b_2)$ form a 312 pattern in $D$ as shown in Figure 12.

$$\underbrace{\begin{pmatrix} 0 & 0 & 0 & \boxed{1} & 0 \\ 0 & 1 & 0 & 0 & 0 \\ 1 & -1 & 0 & 0 & 1 \\ 0 & \boxed{1} & 0 & 0 & 0 \\ 0 & 0 & \boxed{1} & 0 & 0 \end{pmatrix}}_{C} \to \underbrace{\begin{pmatrix} 0 & 0 & 0 & \boxed{1} & 0 \\ 0 & 1 & 0 & 0 & 0 \\ 0 & 0 & 0 & 0 & 1 \\ \boxed{1} & 0 & 0 & 0 & 0 \\ 0 & 0 & \boxed{1} & 0 & 0 \end{pmatrix}}_{D}$$

**Fig. 12:** When only the smallest entry of a classical 312 pattern (formed by the removal of a 321) is affected by the removal of a $-1$ in the SW key process, a 312 is formed in the same rows as the original pattern.

This proves that if a 312 pattern has been created at some stage of the SW key process from $A$ that classically contains 321, then it cannot be removed at a later stage, which concludes the proof of Lemma 4.5. $\square$

Combining the results of Lemma 4.3 and Lemma 4.5, we have the following result.

**Proposition 4.6.** *Let $A$ be an ASM that key-avoids* $312$*. Then $A$ also key-avoids* $321$ *if and only if $A$ classically avoids* $321$.

Building on Proposition 4.6, we can identify the monotone triangles in the image of $\mathrm{Av}_n^{\mathrm{key}}(312, 321)$ through the bijection of Section 3.

**Theorem 4.7.** *The class of $n \times n$ ASMs whose key avoids both* $321$ *and* $312$ *is in bijection with gapless monotone triangles of size $n$ with at most two distinct values in any column.*

**Proof:** Recall from Theorem 3.8 that the monotone triangles corresponding to 312-key-avoiding ASMs are those that are gapless. Consider such a gapless monotone triangle $t_A$ corresponding to the ASM $A$ whose key avoids both 321 and 312.

Suppose $t_A$ has three distinct values in the $k$th column. Since the smallest entry of the column must be $k$ as the columns decrease (gaplessly), there must be an entry equal to $k+2$ in the $k$th column as shown in Figure 13. The row $r$ containing $k+2$ in the $k$th column is missing two smaller values from $\{1, 2, 3, \ldots, k+1\}$, say $a, b$ where $a < b$. Translating back to $A$, this implies that there are 1s in columns $a, b$ that are lower than at least one 1 in column $k+2$. Thus there exist $r < x \leq y$ such that the $A$ has 1 entries in positions $(r, k+2), (x, b), (y, a)$ or such that the $A$ has 1 entries in positions $(r, k{+}2), (x, a), (y, b)$. Choose $k, a, b$ such that $(r, x, y)$ is the largest triple that meets all of the conditions above for the lexicographic order.

$$t_A = \begin{array}{ccccc} \cdot & & & & \\ \cdot & \cdot & & & \\ \vdots & \vdots & \cdot & & \\ \cdot & \cdots & k+2 & \cdots & \\ \vdots & \vdots & \vdots & \cdots & \\ \cdot & \cdots & k+1 & \cdots & \\ \vdots & \vdots & \vdots & \cdots & \\ 1 & \cdots & k & \cdots & n \end{array}$$

**Fig. 13:** A monotone triangle with three distinct values in the $k$th column.

We distinguish cases as illustrated in Fig. 14.

$$\begin{array}{c} \\ r \\ x \\ y \end{array}\!\!\begin{array}{c} \begin{array}{ccc} a & b & k+2 \end{array} \\ \begin{pmatrix} & & 1 \\ & 1 & \\ 1 & & \end{pmatrix} \end{array} \qquad \begin{array}{c} \\ r \\ x \\ z \end{array}\!\!\begin{array}{c} \begin{array}{cccc} a & \gamma & b & k+2 \end{array} \\ \begin{pmatrix} & & & 1 \\ 1 & -1 & 1 & \\ & 1 & & \end{pmatrix} \end{array} \qquad \begin{array}{c} \\ r \\ x \\ y \end{array}\!\!\begin{array}{c} \begin{array}{ccc} a & b & k+2 \end{array} \\ \begin{pmatrix} & & 1 \\ 1 & & \\ & 1 & \end{pmatrix} \end{array}$$

**Fig. 14:** From left to right: Cases 1, 2 and 3 in the proof of Theorem 4.7.

Case 1: If $r < x < y$ and $A$ has 1 entries in positions $(r, k+2), (x, b), (y, a)$, then $A$ classically contains 321, which is impossible by Proposition 4.6.

Case 2: If $r < x = y$ and $A$ has 1 entries in positions $(r, k+2), (x, b), (x, a)$, then due to the $-1$ that must appear in position $(x, \gamma)$ for some $a < \gamma < b$, there is also a 1 in position $(z, \gamma)$ for some $z > x$ and $(r, k+2), (x, b), (z, \gamma)$ forms a classical 321 pattern in $A$. Again, this is impossible from Proposition 4.6.

Case 3: Otherwise, it holds that $r < x < y$ and that $A$ has 1 entries in positions $(r, k+2), (x, a), (y, b)$. In this case, let us show that in going from down from row $x-1$ to row $x$ of the monotone triangle $t_A$, we only insert the value $a$ in row $x$. This amounts to showing that the matrix $A$ has no $-1$ entry on row $x$.

First, there cannot be a $-1$ entry in row $x$ to the left of column $a$, since (among other reasons) the 1 below it would create a 321 pattern in $A$ (with the 1s in positions $(x, a)$ and $(r, k+2)$), which is forbidden by Proposition 4.6.

Next, observe that for any column of index $j \leq k+2$, $j \neq a, b$, the partial sum until row $r$ in column $j$ of $A$ is 1 (because $j$ appears on row $r$ in $t_A$). Using this fact, the alternating sign condition on non-zeros entries of an ASM, and the property that $A$ avoids 321 guaranteed by Proposition 4.6, we deduce that all entries of $A$ in positions $(u, v)$ for $r < u \leq x$ and $a < v \leq k+2$ are 0 (except perhaps in $(x, k+2)$ where we have not yet ruled out the possibility of a $-1$). In passing, this implies that all values from $a+1$ until $k+2$ except $b$ appear on row $x-1$ of $t_A$.

So, assuming that there is a $-1$ in row $x$ of $A$, we obtain that it is in some column $v \geq k+2$, and has a 1 to its right, in some column $d > v$. In row $x$ of $t_A$, at least the values $b$ and $v$ are missing so that in column $d$ of row $x > r$, $t_A$ has a value which is at least $d+2$. If it is exactly $d+2$, this contradicts the maximality of $r$. (Indeed, this difference of 2 implies that there are two 1 entries in $A$ to the southwest of the 1 entry in position $(x, d+2)$, just like we observed with the entry $k+2$ in column $k$ of $t_A$ at the beginning of the proof.) If not, as $t_A$ is gapless, there is a row of $t_A$ of index larger than $x$, such that the entry in column $d$ must be $d+2$, and again this contradicts the maximality of $r$.

Going back to the monotone triangle $t_A$, we have proved that $a$ is the only entry that appears in row $x$ but not in row $x-1$. The values $(a+1)\ldots(b-1)(b+1)\ldots(k+2)$ sitting on columns $a$ to $k$ of row $x-1$ of $t_A$ are then shifted one column to the right on row $x$, as shown in Figure 15. In column $b-1$, we therefore find $b-1$ in row $x$ and $b+1$ immediately above $b-1$ in row $x-1$. However, then the monotone triangle $t_A$ is not gapless in column $b-1$, contradicting $A$ being a 312-key-avoiding ASM by Theorem 3.8. $\square$

$$
t_A = \begin{array}{cccccccccccccc}
\vdots & & & & & & & & & & & & & \\
\vdots & & & & & & & & & & & & & \\
\cdots & \cdots & a-1 & a+1 & \cdots & \cdots & b-1 & b+1 & \cdots & \cdots & k+2 & \cdots & & \\
\cdots & \cdots & a-1 & a & a+1 & \cdots & \cdots & b-1 & b+1 & \cdots & \cdots & k+2 & \cdots & \\
\vdots & \vdots & \vdots & \vdots & \vdots & \vdots & \vdots & \vdots & \vdots & \vdots & \vdots & \vdots & \vdots & \\
1 & \cdots & a-1 & a & a+1 & \cdots & \cdots & b-1 & b & \cdots & k & k+1 & \cdots & n
\end{array}
$$

**Fig. 15:** A monotone triangle showing the insertion of $a$ and eventual insertion of $b$.

From the characterization of $\mathrm{Av}_n^{\mathrm{key}}(312, 321)$ given by Theorem 4.7, we can now prove the main result of this section.

**Theorem 4.8.** *For every $n$, $|\mathrm{Av}_n^{\mathrm{key}}(312, 321)| = \frac{1}{n+1}\binom{2n}{n}$.*

**Proof:** By Theorem 4.7, $\mathrm{Av}_n^{\mathrm{key}}(312, 321)$ is in bijection with the set of gapless monotone triangles of size $n$ with at most two distinct values in any column. To prove the statement, we describe below a size-preserving bijection between gapless monotone triangles with at most two distinct values in any column and another family known to be enumerated by the Catalan numbers.

Recall that an *inversion sequence* $e = e_1e_2e_3 \cdots e_n$ is a sequence of nonnegative integers such that $e_i < i$ for all $i$. An inversion sequence $e$ contains the pattern 10 if there are indices $i < j$ such that $e_i > e_j$. Otherwise, it avoids 10, which is equivalent to $e$ being weakly increasing. Martinez and Savage [11, Section 2.14.1, Theorem 27] proved that the number of inversion sequences avoiding 10 of size $n$ is the $n$th Catalan number $C_n = \frac{\binom{2n}{n}}{n+1}$, for every $n$. (Note that this result is also obtainable by subtracting one from each entry of the sequences defined in Catalan Exercise 78 of [18, p. 32]).

Notice that gapless monotone triangles with at most two distinct values in each column are such that the possible values in column $i$ are only $i$ and $i + 1$ for $i < n$. This is because column $i$ contains value $i$ at the bottom of the triangle (which is the only entry in column $n$) and then the monotonicity and gapless requirements force a second value (if it exists) in the same column to be $i + 1$.

Let $T_n$ be the set of gapless monotone triangles with $n$ columns with at most two distinct values in each column. Let $M_n$ be the set of inversion sequences of size $n$ avoiding 10. Define $f : T_n \to M_n$ to be such that the $i$th entry in the resulting inversion sequence is one less than the number of times $n + 1 - i$ appears in the $(n + 1 - i)$th column. Since there are only $i$ terms in the $(n + 1 - i)$th column and at least one of them must be $n + 1 - i$, the resulting sequence is indeed an inversion sequence. Further, because the monotone triangles have strictly increasing rows, for $i < j$, the number of entries $n + 1 - j$ in the $(n+1-j)$th column is larger than or equal to the number of entries $n+1-i$ in the $(n+1-i)$th column. Hence the resulting inversion sequence is weakly increasing, i.e. avoids 10.

The function $f$ is invertible as the inversion sequence tells us for each column $i$, how many times we repeat the initial column value $i$ beyond the bottom row. As the preimage was a gapless monotone triangle with only two possible column values for each column, the remaining entries in the $i$th column must be $i + 1$. □

We show the bijection given in the proof of Theorem 4.8 applied to a monotone gapless triangle with at most two distinct values in each column in Example 4.9.

**Example 4.9.** In the gapless monotone triangle on the left of Fig. 16, notice that there are:

- 0 repetitions of 5 above the bottom row in column 5,
- 0 repetitions of 4 above the bottom row in column 4,
- 1 repetition of 3 above the bottom row in column 3,
- 1 repetition of 2 above the bottom row in column 2, and
- 3 repetition of 1 above the bottom row in column 1,

giving us the inversion sequence 00113.

**Remark 4.10.** Note that given a gapless monotone triangle $M$ with at most two values in each column, we may also easily draw the Dyck path of the corresponding inversion sequence in the bijection of the proof of Theorem 4.8. Draw a horizontal line between value $i$ and $i + 1$ in column $i$ of $M$, noting that there may be no $i + 1$ value, in which case draw a horizontal line above all values. Then connect the horizontal lines with vertical lines between the columns. Finally, add one additional vertical line to the

$$\begin{matrix} 2 & & & & \\ 1 & 3 & & & \\ 1 & 3 & 4 & & \\ 1 & 2 & 3 & 5 & \\ 1 & 2 & 3 & 4 & 5 \end{matrix} \quad \xrightarrow{f} \quad 00113 \qquad ;$$

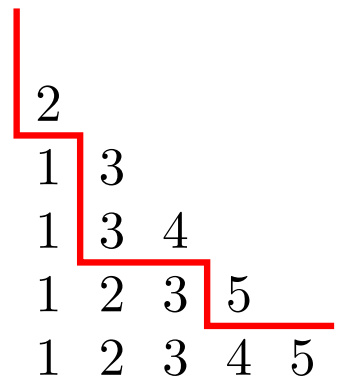

**Fig. 16:** Left: The bijection between gapless monotone triangles with at most two values in each column and inversion sequences avoiding 10, on the example of Example 4.9. Right: The corresponding Dyck path, as described in Remark 4.10.

left of the first column starting at the left side of the first horizontal edge and continuing until one row above the topmost row. This forms a Dyck path (which may be rotated w.r.t. the usual convention). See an example in Fig. 16, right.

# 5 Products of Catalan numbers and enumeration of ASMs with a given key avoiding 312 and 321

In the prior section, we showed that the number of ASMs in $\mathrm{Av}_n^{\mathrm{key}}(312, 321)$ is the $n$th Catalan number $C_n = \frac{1}{n+1}\binom{2n}{n}$. In this section, we show in Theorem 5.5 that the number of ASMs whose SW key is $\pi$ for any given $\pi \in \mathrm{pAv}(312, 321)$ is a product of Catalan numbers. Thus, we give a combinatorial proof of an identity relating a Catalan number to a sum of products of Catalan numbers (see Corollary 5.6). We will need the following definitions and lemma.

**Definition 5.1.** Given matrices $A$ and $B$, define their *direct sum* as the block-diagonal matrix $A \oplus B = \begin{pmatrix} A & 0 \\ 0 & B \end{pmatrix}$ and their *skew-sum* as the block-antidiagonal matrix $A \ominus B = \begin{pmatrix} 0 & A \\ B & 0 \end{pmatrix}$.

**Definition 5.2.** For any $n \geq 1$, let $w_n$ denote the permutation $23 \ldots n1$ (*i.e.*, the skew-sum of the increasing permutation of size $n-1$ and the permutation of size $1$).

**Definition 5.3.** Let $\pi$ be a permutation and $A$ its permutation matrix. We define a *southwest-record* of $\pi$ as a $1$ entry in $A$ that has no $1$ entry to its southwest. (Note that these correspond to right-to-left minima in the one-line notation of the permutation, but we won't use this fact.)

**Lemma 5.4.** *For any permutation* $\pi \in \mathrm{pAv}_n(312, 321)$*, there exists some* $k \geq 1$ *and a composition* $(m_1, m_2, \ldots, m_k)$ *of* $n$ *into* $k$ *(non-zero) parts such that* $\pi = w_{m_1} \oplus w_{m_2} \oplus \cdots \oplus w_{m_k}$*. In other words, the permutation matrix of* $\pi$ *consists of blocks* $w_{m_1}, w_{m_2}, \ldots, w_{m_k}$ *placed along the northwest-southeast diagonal.*

*This correspondence between* $\mathrm{pAv}_n(312, 321)$ *and compositions of* $n$ *(into non-zero parts) is a bijection.*

**Proof:** Let $\pi \in \mathrm{pAv}_n(312, 321)$, and denote by $A$ its permutation matrix. The southwest records of $\pi$ can be ordered from left to right in $A$, as $r_1, r_2, \ldots r_k$ for some $k \geq 1$. Because of the avoidance of 312, for any $j < k$, there is no 1 entry in $A$ that lies both above $r_j$ and to the right of $r_{j+1}$. This forces $\pi$ to be decomposable into a direct sum of $k$ blocks, each containing one southwest-record. The fact that each of the blocks in this direct sum is a permutation of the form $w_m$ for some $m \geq 1$ is a immediate consequence of the avoidance of 321. This is illustrated in Fig. 17, and proves the first statement of the lemma.

The second statement of the lemma follows since for every composition $(m_1, m_2, \ldots, m_k)$ of $n$ into $k$ (non-zero) parts, the permutation $w_{m_1} \oplus w_{m_2} \oplus \cdots \oplus w_{m_k}$ is in $\mathrm{pAv}_n(312, 321)$. $\square$

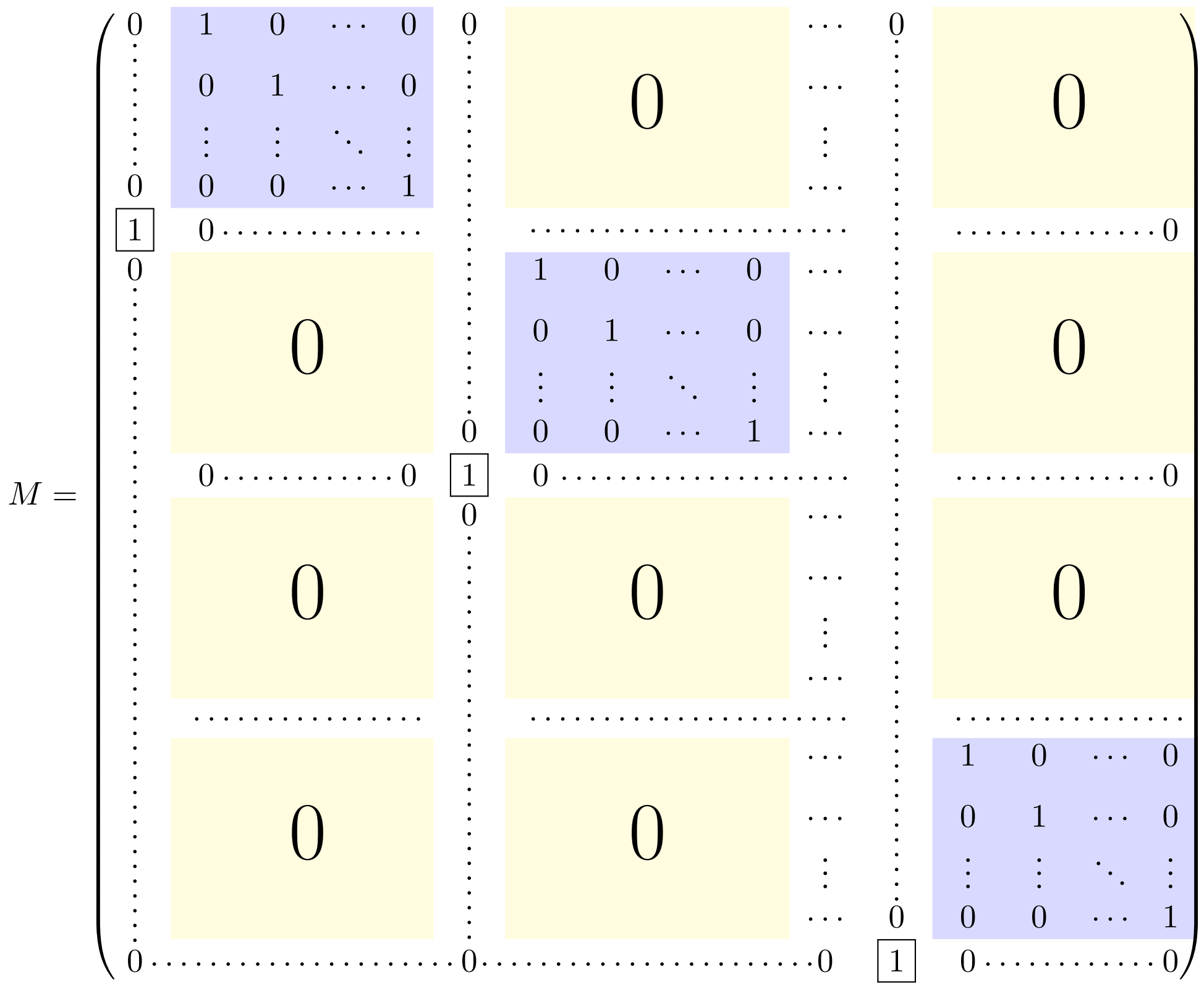

**Fig. 17:** A permutation matrix avoiding 321 and 312 with southwest-records boxed.

We now use Lemma 5.4 along with a result from [20] about the *Schubert polynomial* of $w_n$ to prove our main theorem of this section. (The definition of the Schubert polynomial will not be needed for this proof; we refer the interested reader to the references given.)

**Theorem 5.5.** *Given $\pi \in \mathrm{pAv}(312, 321)$, the number of ASMs whose SW key is $\pi$ is a product of Catalan numbers. In particular, if $\pi = w_{m_1} \oplus w_{m_2} \oplus \cdots \oplus w_{m_k}$, then the number of ASM whose key is $\pi$ equals $\prod_{i=1}^{k} C_{m_i - 1}$.*

**Proof:** Let $\pi \in \mathrm{pAv}_n(312, 321)$. By Lemma 5.4, $\pi$ is a direct sum of permutations $w_i = 23 \cdots i1$. So suppose $\pi = w_{m_1} \oplus w_{m_2} \oplus \cdots \oplus w_{m_k}$.

Consider starting with the permutation $\pi$ and creating an ASM from this key by successively adding negative ones using the reverse key process. The only expansions that are possible starting with the permutation matrix $\pi$ are to take the southwest record (see Definition 5.3) in any block and expand it to a $-1$ by a *simple expansion* (inverse of the *simple removal* from Definition 1.3). It is simple, because any removal of a $-1$ (simple or not) in the SW-key process must keep the overall sum of the matrix the same. Thus, the removal of a $-1$ must also result in one fewer $1$. Therefore, the expansion of a single $1$ must result in a $-1$ with a $1$ to the west and south of it. This is the only possible expansion, because any $-1$ added by an expansion must have a $1$ to the north and east of it already. So none of the 1s other than the southwest record may expand to add a $-1$. For the same reason, this and all subsequent expansions involving these $1$ entries stay within a block of the decomposition $\pi = w_{m_1} \oplus w_{m_2} \oplus \cdots \oplus w_{m_k}$, so the choices of expansions in different blocks are independent.

Thus, the number of ASMs whose key is $\pi$ for $\pi = w_{m_1} \oplus w_{m_2} \oplus \cdots \oplus w_{m_k} \in \mathrm{pAv}(312, 321)$ is given by the product of the number of ASMs in whose key is $w_{m_i}$, $1 \leq i \leq k$. It remains to show that these are given by the $(w_{m_i} - 1)$st Catalan numbers.

By [20, Prop. 1], the *Schubert polynomial* with all variables set equal to 1 of the permutation $w_n$ equals the $(n-1)$st Catalan number. In symbols: $\mathfrak{S}_{w_n}(1, \ldots, 1) = C_{n-1}$. So, it is actually enough to prove that the number of ASMs whose key is $w_n$ is $\mathfrak{S}_{w_n}(1, \ldots, 1)$.

By [10], the Schubert polynomial $\mathfrak{S}_\sigma$ of a permutation $\sigma$ of size $n$ is a sum of positive monomials indexed by objects called *reduced bumpless pipe dreams* on an $n \times n$ grid, where the sum is taken over all of these objects with associated permutation $\sigma$. Bumpless pipe dreams (reduced and non-reduced) on an $n \times n$ grid with associated permutation $\sigma$ are in direct bijection with $n \times n$ alternating sign matrices with key $\sigma$ [19]. So if we show that the bumpless pipe dreams corresponding to ASMs with key $w_n$ are all reduced, we can use the above-mentioned result of [20] to count them.

Now, since $w_n$ avoids 312, $w_n$ also avoids 3412. It was shown in [19, Lemma 7.2] that the bumpless pipe dreams for a permutation $\sigma$ are all reduced if $\sigma$ avoids 3412. (Note that [19, Lemma 7.2] uses a convention equivalent to the northwest key, and thus states the pattern as 2143. Since we use the southwest key, the avoided permutation for us is the reverse.) Thus, the bumpless pipe dreams corresponding to ASMs with key $w_n$ are reduced, so there are as many ASMs with key $w_n$ as there are terms in the Schubert polynomial $\mathfrak{S}_{w_n}$, concluding the proof. □

It can be noted that the arguments of the above proof also show that there are $\mathfrak{S}_\sigma(1, \ldots, 1)$ ASMs with key $\sigma$ for any permutation $\sigma$ that avoids 3412.

We introduce the following notation to assist in completing the proof of the Catalan Identity given in Corollary 5.6. Given $\pi \in \mathrm{pAv}_n(312, 321)$, let us denote by $B(\pi)$ the composition $(m_1, m_2, \ldots, m_k)$ identified in Lemma 5.4. Furthermore, let $\mathfrak{C}_{n,k}$ denote the set of weak compositions of $n$ into $k$ parts (meaning that we allow parts that are equal to 0), and let $\mathfrak{C}^{>0}_{n,k}$ denote the set of strict compositions of $n$ into $k$ parts (meaning that all parts are greater than 0). For example,

$$\mathfrak{C}_{4,2} = \{(4,0), (3,1), (2,2), (1,3), (0,4)\}, \text{ while } \mathfrak{C}^{>0}_{4,2} = \{(3,1), (2,2), (1,3)\}$$

If $m \in \mathfrak{C}_{n,k}$ or $\mathfrak{C}^{>0}_{n,k}$, we write $m = (m_1, m_2, \ldots, m_k)$.

**Corollary 5.6.**

$$C_n = \sum_{\pi \in \mathrm{pAv}_n(312,321)} \prod_{i \in B(\pi)} C_{i-1} \tag{1}$$

$$= \sum_{k=1}^{n} \sum_{m \in \mathfrak{C}_{n-k,k}} C_{m_1} C_{m_2} \cdots C_{m_k} \tag{2}$$

**Proof:** The first part of this corollary follows directly from Theorems 4.8 and 5.5. Theorem 4.8 says the number of $n \times n$ ASMs whose key avoids 312 and 321 equals $C_n$, while Theorem 5.5 says the number of $n \times n$ ASMs with key $\pi$ where $\pi$ avoids 312 and 321 equals $\prod_{i \in B(\pi)} C_{i-1}$. So summing the latter over all permutations in $\mathrm{pAv}_n(312, 321)$ equals $C_n$.

The second part of the corollary comes from the second part of Lemma 5.4, which says each strict composition of $n$ determines a permutation in $\mathrm{pAv}_n(312, 321)$. So instead of summing over permutations avoiding 312 and 321, we sum over strict compositions of $n$, obtaining:

$$C_n = \sum_{k=1}^{n} \sum_{t \in \mathfrak{C}^{>0}_{n,k}} C_{t_1-1} C_{t_2-1} \cdots C_{t_k-1}.$$

If we subtract one from each part of $t \in \mathfrak{C}^{>0}_{n,k}$, we obtain a weak composition of $n-k$ into $k$ parts. This is a bijection between $\mathfrak{C}^{>0}_{n,k}$ and $\mathfrak{C}_{n-k,k}$. Thus, by shifting the indices on each Catalan number, we can sum over weak compositions and obtain (2). □

We note that the equality (2) can also be seen as a consequence of the recurrence

$$C_n = \sum_{i=1}^{n} C_{i-1} C_{n-i}.$$

More precisely, (2) follows by induction from the above recurrence. Checking the base case for $n = 1$ is obvious. Then, for $n > 1$, we can write

$$\begin{aligned}
C_n &= \sum_{i=1}^{n} C_{i-1} C_{n-i} \\
&= C_{n-1} + \sum_{i=1}^{n-1} C_{i-1} \left( \sum_{j=1}^{n-i} \sum_{m \in \mathfrak{C}_{n-i-j,j}} C_{m_1} C_{m_2} \cdots C_{m_j} \right) \quad \text{by the induction hypothesis} \\
&= C_{n-1} + \sum_{i=1}^{n-1} C_{i-1} \left( \sum_{k=2}^{n-i+1} \sum_{m \in \mathfrak{C}_{n-i-k+1,k-1}} C_{m_1} C_{m_2} \cdots C_{m_{k-1}} \right) \\
&\qquad \text{with the change of variables } k = j+1
\end{aligned}$$

$$= C_{n-1} + \sum_{k=2}^{n} \sum_{m \in \mathfrak{C}_{n-k,k}} C_{m_1} C_{m_2} \cdots C_{m_k} \quad \text{by setting } m_k = i-1$$

$$= \sum_{k=1}^{n} \sum_{m \in \mathfrak{C}_{n-k,k}} C_{m_1} C_{m_2} \cdots C_{m_k},$$

thus proving (2).

## Acknowledgments

This work began at the 2023 Schloss Dagstuhl workshop on Pattern Avoidance, Statistical Mechanics and Computational Complexity. We wish to thank the institute and the workshop organizers for the conducive research environment and the referees for helpful comments. Smith thanks the Mocqua team at Loria for supporting her visit in Nancy.